\documentclass[twoside,twocolumn]{article}

\usepackage{blindtext} 

\usepackage[sc]{mathpazo} 
\usepackage[T1]{fontenc} 
\usepackage{microtype} 

\usepackage[english]{babel} 

\usepackage[hmarginratio=1:1,top=32mm,columnsep=20pt]{geometry} 
\usepackage[hang, small,labelfont=bf,up,textfont=it,up]{caption} 
\usepackage{booktabs} 

\usepackage{lettrine} 

\usepackage{enumitem} 
\setlist[itemize]{noitemsep} 

\usepackage{abstract} 

\usepackage{titlesec} 
\renewcommand\thesection{\Roman{section}} 
\renewcommand\thesubsection{\roman{subsection}} 
\titleformat{\section}[block]{\large\scshape\centering}{\thesection.}{1em}{} 
\titleformat{\subsection}[block]{\large}{\thesubsection.}{1em}{} 

\usepackage{fancyhdr} 
\usepackage{titling} 

\usepackage{hyperref} 

\pretitle{\begin{center}\Huge\bfseries} 
\posttitle{\end{center}} 
\title{Necessary and sufficient conditions for exponential stabilization of  systems  affine in controls} 
\author{%
\textsc{M. Ouzahra}\\
\normalsize MASI Team, University of
Sidi Mohamed Ben Abdellah \\ 
\normalsize \href{mohamed.ouzahra@usmba.ac.ma}{mohamed.ouzahra@usmba.ac.ma} 
}
\date{}
\renewcommand{\maketitlehookd}{%
\begin{abstract}
This paper provides  necessary and sufficient conditions for exponential stabilization of distributed
systems affine in control, evolving in a Banach state space, by means of constant controls. An explicit  estimate of the convergence speed is obtained under the given constant  control. Various applications are provided.
\end{abstract}
}

\newtheorem{thm}{Theorem}
\newtheorem{lem}[thm]{Lemma}
\newtheorem{rem}[thm]{Remark}
\newtheorem{ex}[thm]{Example}

\newtheorem{cor}[thm]{Corollary}
\newtheorem{df}[thm]{Definition}
\newtheorem{pf}[thm]{Proof}

\begin{document}

\maketitle


\section{Introduction}

In this work, we deal with the following control-affine  system:
\begin{equation}\label{SS}
\displaystyle\frac{d z(t)}{d t}= Az(t) + v(t)Bz(t),\;\;z(0) = z_0,
\end{equation}
where $v(t)$ is a real valued control, $A$ is the infinitesimal generator of a linear $C_0-$semigroup of contractions
 $S_0(t)$ on a real Hilbert space $H$ with inner product and corresponding norm denoted respectively by $\left<\cdot,\cdot\right > $
and $\|\cdot\|$, so that $A$ is dissipative, i.e. $\langle Az,z\rangle \le 0,$  for all $z\in  {\cal D}(A). $ Here $ B $ is a (possibly) nonlinear
operator from $H$ to $H$ such that $B(0) = 0, $ so that $0$ is an equilibrium
for (\ref{SS}). An important special case of (\ref{SS}) is when  $B$ is a bounded linear operator.\\
 There are numerous real-world problems that can be represented by the system (\ref{SS}). They include  applications in nuclear, thermal, chemical, social processes, etc.. (see \cite{bal1,beau,kha98,khap,moh73}). Feedback stabilization of systems affine in control has been investigated by numerous authors using various control approaches, such as quadratic control laws, sliding mode control, piecewise constant feedback and optimal control laws (see \cite{bacc,bank,bal1,beau,ber10,jurj,ouz10,ouz12a,sle78}). The most popular feedback control for stabilization problem of system (\ref{SS}) is given by:
\begin{equation}\label{quad}
    v(t)=-\langle Bz(t),z(t)\rangle
\end{equation}
In  \cite{bal1}, it has been shown that under
the condition:
\begin{equation}\label{bal}
\langle BS_0(t)y,S_0(t)y\rangle =0,\ \forall t\geq0\Longrightarrow
y=0,
\end{equation}
the quadratic feedback (\ref{quad}) weakly stabilizes the system (\ref{SS}) provided that  $B$ is sequentially continuous from $H_{w}$ ($H$ endowed
with the weak topology) to $H$. Moreover, under the assumption
\begin{equation}\label{coer}
\int_{0}^{T}|\langle BS_0(t)y,S_0(t)y\rangle|dt\geq    \delta\|y\|^{2}, \ \forall y\in H,
\end{equation}
for some constant $ T,\delta>0, $ a strong stabilization result has been obtained using the control (\ref{quad}) (see \cite{ber99,ouz08}). However,
in this way the convergence of the resulting closed loop state is not better than
$\|z(t)\|=\mathcal{O}(\displaystyle\frac{1}{\sqrt{t}}).$ The problem of exponential stabilization of the system  (\ref{SS}) has been considered in \cite{che,ouz10,ouz11m} with the following   bounded  feedback
\begin{equation} \label{e(t)}
v(t)=-\rho \displaystyle\frac{\langle z(t),Bz(t)\rangle }{\|z(t)\|^2}{\bf 1}_{\{t\ge 0;\, z(t)\ne 0\}},
\end{equation}
where $\rho>0$ is the gain control. Moreover, under the assumption (\ref{coer}), the exponential stabilization of (\ref{SS}) has been studied in \cite{ouz12a} using the switching control
\begin{equation} \label{sw(t)}
w(t)=-\rho \:sign(\left<z(t),Bz(t)\right>),\; \rho>0.
\end{equation}
In the case of parabolic like bilinear systems, one may investigate the relation between the
stability of a distributed parameter system and that of a
finite-dimensional one. This idea has been used via a decomposition of the state space according to spectral properties
 of the considered system (see \cite{ouz09,ouz11}).\\
 Among control laws that present more advantages in theory and application, we mention  constant controls (see \cite{bacc2,bacc,kal,liu,rob,tsin}). Indeed, such a control is simple to implement, since it does not require the knowledge of the system's state. Another important point is that this control does not depend on the initial state, so it can be applied to stabilization problem with a priori constraint on the  control.  In the finite dimensional case,  various necessary and sufficient conditions for the stabilization of  system (\ref{SS}) by constant
controls have been formulated in terms of Lyapunov functions (see e.g., \cite{bacc2,bacc,kal,tsin}). In \cite{bacc} it has been proved that every constant  stabilizable bilinear system (\ref{SS}) admits a quadratic Lyapunov function. Conversely, it has  been shown in \cite{tsin} that system (\ref{SS}) is constant stabilizable if a quadratic control Lyapunov function $V(z) = z^TPz $ exists and all the eigenvalues of the symmetric matrix
 $PB + B^TP$ have the same sign. In \cite{kal}, it has been showed that system (\ref{SS}) is constant  stabilizable if  the real parts of the
eigenvalues of $B$ have all the same sign. In \cite{bacc2}, eigenvalues of $B$ are allowed to have real parts
with opposite sign, provided that  some local estimations of extremum  of the quadratic forms $q_1(z) = z^TPAz$ and $q_2(z) = z^TPBz $ are available. In \cite{rob}, the problem of constant stabilization has been considered by using the analysis of the Lie algebra. In the context of infinite-dimensional control-affine system, the authors in \cite{liu} gave sufficient conditions for constant exponential stabilization when $A$ is skew adjoint and $B$ is a bounded linear operator. They also formulated  necessaries conditions when, in addition, $B$ is self adjoint  and dissipative.\\
The aim of this paper is to provide necessary and sufficient conditions for constant exponential stabilizability of finite and infinite dimensional systems that can be described by the system (\ref{SS}). \\The paper is organized as follows : In the second section,  we provide sufficient conditions for exponential stabilization of system (\ref{SS})
with constant controls. In the third section, we  give necessary and sufficient conditions for uniform exponential stabilization by means of  constant controls. In the fourth section, we examine the finite dimensional  case, and we give applications to some class of infinite dimensional systems. The question of robustness is discussed in
the fifth section. Finally, the sixth section is devoted to some applications.



\section{Sufficient conditions for exponential stabilization}

Let us  recall  the following definition of  exponential stabilization of
 system (\ref{SS}) \cite{bacc,zab95}:

\begin{df}
The system (\ref{SS}) is  exponentially stabilizable if there exists a
 feedback control $v(t)= f(z(t)),\; f : H\to R$ such that system (\ref{SS})  satisfies the following
property :\\
(i) for each $z_0$ there exists a unique mild solution $z(t),$ defined for all
 $t\in \mathbf{R}^+, $ of system (\ref{SS}),\\
(ii) there exist $M, \sigma>0$ (depending, eventually, on $z_0$) such that the mild solution $z(t)$ starting at $z_0$ satisfies

\begin{equation}\label{exp}
\|z(t)\| \le M e^{-\sigma t}\|z_0\|,\; \forall t\ge0\cdot
\end{equation}
The system (\ref{SS}) is uniformly  exponentially stabilizable if
(\ref{exp}) holds for some $M$ and $\sigma$ which are independent
of $z_0.$\\
If the closed loop operator  $f$ is a constant map, then we say that system (\ref{SS}) is constant (uniformly) exponentially stabilizable.\\
\end{df}

\begin{rem}
Notice that, unlike the conventional controls (\ref{quad}), (\ref{e(t)}) and (\ref{sw(t)}), the assumption (\ref{coer})  is not sufficient for
 constant  stabilization of (\ref{SS}), as we can see by taking:  $A=0$ and $Bz=|z|,$ in $H:=R.$
\end{rem}
According to the above remark,  it is natural to impose a stronger assumption than (\ref{coer}) in order to get constant
 exponential stabilization. The following theorem provides sufficient conditions for exponential stabilization  with constant controls.
\begin{thm}\label{thm-suff}
Let $A$ generate a linear $C_0-$semigroup $S_0(t)$ of contractions on $H,$  let  $ B : H\rightarrow H$ be positive (i.e. $\langle Bz,z\rangle \ge 0,
\; \forall z\in H$) and lipschitz on any bounded set
 of $H$, and assume that (\ref{coer}) holds.\\ Then for any $R>0,$  there exists $ \lambda_{R}>0$ such that for all $0<\lambda <\lambda_{R},$ the control
\begin{equation}\label{cte0}
 v(t)=- \lambda
\end{equation}
allows the estimate $\|z(t)\|\le M_Re^{-\sigma_R t} \|z_0\|,$ for any $z_0\in H$ with $ \|z_0\|\le R,$   (for some $M_R, \sigma_R >0$).\\
\end{thm}

\begin{pf} The mild solution of (\ref{SS}), whose local existence is guaranteed by the local Lipschitz of $B, $ is given by the following variation of constants formula :
\begin{equation}\label{mild}
  z(t) = S_0(t)z_0 -\lambda \displaystyle\int_0^tS_0(t-s)Bz(s)ds,
\end{equation}
and since $\langle Bz,z\rangle \ge 0,\; $ for all $z\in H, $ the solution $z(t)$ is global (see \cite{Pazy}, p. 185). Furthermore, the map $z_0 \mapsto S_\lambda(t)z_0 :=z(t)$
defines a nonlinear semigroup on $H$.\\
Since $S_0(t)$ is a contraction semigroup, we get by proceeding as in (\cite{bal2}, Lemma 5.5)
\begin{equation}\label{integral}
\begin{array}{cc}
   \|S_\lambda(t)z_0\|^2  - \|S_\lambda(s)z_0\|^2 \le &  \\
    -
2\lambda \int_s^t \langle BS_\lambda(\tau)z_0,S_\lambda(\tau)z_0\rangle d\tau, &  \forall 0\le s\le t.
  \end{array}
\end{equation}
This implies that:
\begin{equation}\label{born}
\|S_\lambda(t)z_0\| \le \|z_0\|,\, \forall z_0\in H.
\end{equation}
From (\ref{mild}) and   (\ref{born}), it comes
\begin{equation}\label{inq1}
\|S_\lambda(t)z_0-S_0(t)z_0\|\le T L_{R} \lambda \|z_0\|,\; \forall t\in [0,T],
\end{equation}
where $L_R$ is a Lipschitz constant of $B$ on the ball $b(O,R)$ of center $O$ and radius $R>0$.\\
Moreover, we can easily verify that for all $t\ge 0,$
\begin{equation}\label{**}
\begin{array}{cc}
  \left<BS_0(t)z_0,S_0(t)z_0\right> = &  \left<BS_\lambda(t)z_0,S_\lambda(t)z_0\right>+ \\
   & \left<BS_0(t)z_0-
BS_\lambda(t)z_0,S_0(t)z_0\right> + \\
   &  \left<BS_\lambda(t)z_0,S_0(t)z_0-S_\lambda(t)z_0\right>.
\end{array}
 \end{equation}
Based on this expression and using (\ref{inq1}), we obtain
$$
\left<BS_0(t)z_0,S_0(t)z_0\right>  \le \left<BS_\lambda(t)z_0,S_\lambda(t)z_0\right> + 2\lambda T L^2_{R} \|z_0\|^2,
$$
which gives after integration
$$
\begin{array}{cc}
  \int_0^T \langle BS_0(t)z_0,S_0(t)z_0\rangle dt \le& 2\lambda T^2 L^2_{R} \|z_0\|^2+ \\
   &   \int_0^T \langle BS_\lambda(t)z_0,S_\lambda(t)z_0\rangle dt.
\end{array}
 $$
Taking $S_\lambda(t)z_0$ instead of $z_0$ in this last inequality, we can see   by using (\ref{coer})  and  the superposition  property of the semigroup $S_\lambda(t)$ that for all $ t\ge 0,$ we have:
\begin{equation}\label{coeLambda1}
\displaystyle\int_t^{t+T} \left<BS_\lambda(s)z_0,S_\lambda(s)z_0\right>  ds \ge (\delta -2\lambda T^2 L^2_{R})\|S_\lambda(t)z_0\|^2.
\end{equation}
In the sequel, we take $ 0<\lambda < \frac{\delta}{2T^2 L^2_{R}}$ and for all $k\in I\!\!N,$ we set $  U_k=\|S_\lambda(kT)z_0\|^2$.\\
It follows from  (\ref{integral}) that for all $k\in I\!\!N, $ we have:
$$U_{k+1} - U_k\le
-2\lambda\displaystyle\int_{kT}^{(k+1)T} \langle BS_\lambda(t)z_0,S_\lambda(t)z_0\rangle dt.
$$
This inequality, together with (\ref{coeLambda1}), gives  :
\begin{equation}\label{Uk}
U_{k+1} \le  \gamma U_{k},
\end{equation}
where $ \gamma = 1-2\lambda(\delta-2\lambda T^2 L_{R}^2).$ Let $\eta>0$ such that $ \gamma\in (0,1),\; $ for all $ \lambda \in (0,\eta).$\\
Thus for
$0<\lambda<\lambda_R:=\inf(\eta,\frac{\delta}{2T^2L^2_{R}}),$ we have :
$$ U_k\leq \gamma^k \|z_0\|^2, \;\forall \; k\ge
0.$$
Let $t\ge0$ and let $k=E(\displaystyle\frac{t}{T})$. We have
$$
 \|S_\lambda(t)z_0\|^2 \le \|S_\lambda(t-kT)\|^2 U_k$$
 $$
 \le \gamma^k \|z_0\|^2.
 $$
It follows that $$
  \|S_\lambda(t)z_0\|\le M_Re^{-\sigma_R t}\|z_0\|, \;\forall t\ge 0,\; \forall z_0\in H,
$$
where  $    \sigma_R = \frac{-\ln \gamma}{2T} $ and $M_R=\frac{1}{\sqrt \gamma},$ which depend, via $\gamma$, on $\|z_0\|$.
\end{pf}

\begin{rem} \begin{enumerate}
\item The result of the above theorem remains true when $B$ is  dissipative, provided that the control (\ref{cte0}) is  replaced by
$v(t)=\lambda>0.$

\item Note that if $B$ is linear, then  the  estimate of Theorem \ref{thm-suff} holds for all initial state $z_0\in H$.
\end{enumerate}
\end{rem}

\section{Uniform exponential stabilization}

\subsection{Sufficient conditions for uniform exponential stabilization}

The following theorem provides sufficient conditions for uniform exponential stabilization  with constant controls.
\begin{thm}\label{thm-suffunif}
Let $A$ generate a linear $C_0-$semigroup $S_0(t)$ of contractions on $H,$ and let $ B$ be a Lipschitz operator such that:
\begin{equation}\label{coerexp}
 \displaystyle\int_0^T \langle B S_0(t)y,S_0(t)y\rangle dt \ge \delta \|y\|^2,
\;\forall y \in H,
\end{equation}
for some $ \delta, T>0.$ Then there exists $ \lambda_{max}>0$, which is independent of the initial state $z_0,$ such that for all   $\lambda \in (0, \lambda_{max}); $   the control
\begin{equation}\label{cte}
 v(t)=- \lambda
\end{equation}
 uniformly  exponentially stabilizes (\ref{SS}), i.e., there exist $M, \sigma >0$ such that
\begin{equation}\label{exp-est}
\|z(t)\|\le M e^{-\sigma t} \|z_0\|,\, \forall z_0\in H.
\end{equation}
\end{thm}

\begin{pf}
Since $B$ is Lipschitz, we have that for any $\lambda>0$ the system (\ref{SS}) admits a unique global mild solution $S_\lambda(t)z_0$. Using Gronwall inequality and the fact that $S_0(t)$ is a contraction semigroup, we deduce from (\ref{mild}) that
\begin{equation}\label{Slambda}
\|S_\lambda(t)z_0\|\le e^{\lambda L t} \|z_0\|,\; \forall z_0\in H,\; \forall t\ge 0,
\end{equation}
where  $L$ denotes a Lipschitz constant of $B$.\\
It follows from (\ref{mild}) and (\ref{Slambda}) that for all $t\in [0,T],$
$$
\|S_0(t)z_0-S_\lambda(t)z_0\|\le (e^{\lambda TL}-1)  \|z_0\|.
$$
Since $e^{\lambda TL}-1\sim \lambda TL $ as $\lambda\to 0,$  there exist $\eta_1, K>0 $ which are independent of $z_0$, such that for all $\lambda\in (0,\eta_1), $ we have\\
\begin{equation}\label{inq}
\|S_0(t)z_0-S_\lambda(t)z_0\|\le K \lambda \|z_0\|,\; \forall t\in [0,T].
\end{equation}
Using  (\ref{Slambda}) and  (\ref{inq}), we obtain from (\ref{**})
$$
\left<BS_0(t)z_0,S_0(t)z_0\right>  \le \left<BS_\lambda(t)z_0,S_\lambda(t)z_0\right> + \lambda \tilde{K}\|z_0\|^2,
$$
for all $ t\in [0,T]$ and $\lambda\in (0,\eta_1), $ where $\tilde{K}=2KL(1+e^{\eta_1L T}).$\\
This gives after integration
$$
\begin{array}{cc}
  \int_0^T \langle BS_0(t)z_0,S_0(t)z_0\rangle dt \le & \lambda T \tilde{K} \|z_0\|^2 + \\
   &   \int_0^T \langle BS_\lambda(t)z_0,S_\lambda(t)z_0\rangle dt.
\end{array}
$$
It follows from this last inequality and (\ref{coerexp}) that  for all $ t\ge 0,$
\begin{equation}\label{coeLambda}
\displaystyle\int_t^{t+T} \left<BS_\lambda(s)z_0,S_\lambda(s)z_0\right>  ds \ge (\delta -\lambda T \tilde{K})\|S_\lambda(t)z_0\|^2.
\end{equation}
Then, the remaining part of the proof is similar to the one performed for Theorem \ref{thm-suff}. As a consequence we obtain the estimate
(\ref{exp-est})  for
$  \sigma = \frac{-\ln \gamma}{2T} $ and $M=\frac{1}{\sqrt \gamma} \displaystyle\sup_{0\le t\le T}\|S_\lambda(t)\| $
 (which are independent of $z_0$), provided that
$0<\lambda< \lambda_{max} : = \inf{(\eta_1,\eta_2,\frac{\delta}{T \tilde{K}})}, $ where $\eta_2$ is such that $\gamma := 1-2\lambda(\delta-\lambda T \tilde{K})\in (0,1).$
\end{pf}

\begin{rem}\label{e-sw}

\begin{enumerate}

\item If the system (\ref{SS}) is subject to the control constraint $
 |v(t)|\le v_{\max}, $ then one may choose the constant control (\ref{cte}) such that $\lambda \in (0,\min(\lambda_{\max},v_{\max}))$.

    \item It is easily verified, from the expression of $\sigma, $ that the best value of the rate of exponential convergence $\sigma$ corresponds
to $v(t)=-\displaystyle\frac{\delta}{2T\tilde{K}}.$

\item Unlike the non constant controls, the weak assumption (\ref{bal}) is not sufficient for weak stabilization, as evidenced with   the example: $
A=0, \;\; By=|y|, \; H=I\!\!R.$

\item   The assumption that $S_0(t)$ is a contraction semigroup is essential in the statement of Theorem \ref{thm-suffunif}. Indeed, let us consider the system (\ref{SS}) on $H:=R^2$ with $A=\left(
         \begin{array}{cc}
           0 & 1 \\
           0 & 0 \\
         \end{array}
       \right)\; \mbox{and} \; B=\left(
                        \begin{array}{cc}
                          1 & 0 \\
                          0 & -1 \\
                        \end{array}
                      \right).$
We have  $e^{tA}=\left(
         \begin{array}{cc}
           1 & t \\
           0 & 1 \\
         \end{array}
       \right).$ Then for all $z=(x,y)\in H:=\mathbb{R}^2,$ we have  $
\langle Be^{tA}z,e^{tA}z\rangle=x^2+2txy+t^2y^2-y^2.$ Thus the assumption (\ref{coerexp}) holds for some $T$ large enough and $0<\delta< T$. Indeed, it suffices to look for $T$ and $\delta$ such that for all $x, y \in I\!\!R$ we have: $
(T-\delta)x^2+T^2 xy +(\frac{T^3}{3}-T-\delta)y^2 \ge 0.$ For $y=0,$ this inequality is verified provided that $T-\delta\ge0$. In the case $y\ne 0$, we should have $T-\delta>0$  and $\Delta=(T^4-4(T-\delta)(\frac{T^3}{3}-T-\delta)) y^2<0$, which is guaranteed by letting  $T\to +\infty$, since $\Delta \sim -\frac{T^4}{3},$ as $T\to +\infty$. However, the system is not constant stabilizable,
since for all $\lambda\in R,$ we have:$
A-\lambda B=
\left(
                        \begin{array}{cc}
                          -\lambda & 1 \\
                          0 & \lambda \\
                        \end{array}
                      \right).$

       \end{enumerate}

\end{rem}

\subsection{Necessary conditions for  uniform exponential stabilization}

In the following theorem, we  give a necessary condition for  uniform exponential stabilizability of (\ref{SS}) with  constant controls.

\begin{thm}\label{thm-necaffunif}
Suppose that $A$ generates a semigroup of isometries $S_0(t)$ and let the operator $B$ be Lipschitz.\\
If the system (\ref{SS}) is uniformly exponentially stabilizable  with a constant control, then there exist   $T, \;\delta> 0$ such that
\begin{equation}\label{coerexpaff}
 \displaystyle\int_0^T \|B S_0(t)y\| dt \ge \delta \|y\|,
\;\forall y \in H.
\end{equation}
\end{thm}

\begin{pf}
Let $v(t)=-\lambda, \; \lambda\in I\!\!R^+ $ be an uniformly exponentially stabilizing control for (\ref{SS}), and let  $M\ge 1,\, \sigma>0 $
be such that the corresponding mild solution  $ z(t)=S_\lambda(t)z_0$ of (\ref{SS}) satisfies the estimate (\ref{exp-est}). Since $S_0(t)$ is  of isometries, we have  $\langle Az,z\rangle =0, $ for all $z\in {\cal D}(A)$. Then, using a density argument, we show as for (\ref{integral}) that  for all $T>0$ and  all $z_0\in H$ we have:
\begin{equation}\label{estimcn}
 2\lambda  \int_0^{T}\langle B S_\lambda(t)z_0,S_\lambda(t)z_0\rangle dt
= \|z_0\|^2 - \|S_\lambda(T)z_0\|^2.
\end{equation}
Now, remarking that for $z_0=0,$ we have $z(t)=0,\; \forall t\ge 0;$ we can suppose, in the remainder of the proof, that $z_0\ne0.$ Using the  fact that $ S_\lambda(t)z_0 \to 0, $ as $t\to+\infty$, it comes from (\ref{estimcn})  that $\lambda\ne0$.\\
The estimate  (\ref{exp-est}) together with (\ref{estimcn}), gives
\begin{equation}\label{nc0}
 |\int_0^{T}\langle B S_\lambda(t)z_0,S_\lambda(t)z_0\rangle dt|
\ge \frac{1-Me^{-\sigma T}}{2|\lambda|}\|z_0\|^2,\, \forall z_0\in H,
\end{equation}
for  $T>\displaystyle\frac{\ln(M)}{\sigma}. $  From the variation of constant formula, we have
$$
\|BS_\lambda(t)z_0\|\le \|BS_0(t)z_0\| +|\lambda| \|B\|\displaystyle\int_0^t\|BS_\lambda(s)z_0\|ds.
$$
Then the  Gronwall inequality \cite{cor} yields, for $t\in [0,T],$
$$
\|BS_\lambda(t)z_0\|\le \|BS_0(t)z_0\| +|\lambda| \|B\|\displaystyle\int_0^t \|BS_0(s)z_0\|e^{|\lambda| (t-s)\|B\|}ds
$$
$$
\le \|BS_0(t)z_0\| +|\lambda| \|B\| e^{|\lambda| \|B\|T} \displaystyle\int_0^{T} \|BS_0(s)z_0\|ds.
$$
It follows that
$$
\begin{array}{cc}
  \int_0^{T} \|BS_\lambda(t)z_0\| dt \le & \int_0^{T} \|BS_0(t)z_0\| dt + \\
   & |\lambda| \|B\| Te^{|\lambda| T \|B\|} \displaystyle\int_0^{T}  \|BS_0(t)z_0\|dt.
\end{array}
   $$
This inequality, together with  (\ref{nc0}),  implies $$ \displaystyle\int_0^{T} \|B S_0(t)y\| dt \ge \delta \|y\|,
\;\forall y \in H.
$$
for $T>\displaystyle\frac{\ln(M)}{\sigma} $ and $\delta=\frac{1-Me^{-\sigma T}}{2M|\lambda|} \big (1+|\lambda| T \|B\| e^{|\lambda| T \|B\|} \big )^{-1}.$

\end{pf}

\begin{rem}
\begin{enumerate}

\item The inequality (\ref{coerexpaff}) is  not sufficient for constant uniform exponential stabilization, as evidenced  by the example :
$    A=0 \; \mbox{ and}\;  Bz=|z|,\;  \mbox{in} \; H:=R.$

\item  We have : (\ref{coerexp}) $\Rightarrow$ (\ref{coer}) $\Rightarrow$ (\ref{coerexpaff}).

\item In general, even if $\dim H<\infty, $ we have

$\bullet$  (\ref{coer}) $\not\Rightarrow$ (\ref{coerexp}), as we can see for $A=0$ and $Bz=|z|,$ on $H:=R.$

$\bullet$  (\ref{coerexpaff}) $\not\Rightarrow$ (\ref{coer}), as we can see by taking : $A=0$ and $ B=\left(
                                                                                       \begin{array}{cc}
                                                                                         0 & 1 \\
                                                                                         -1 & 0 \\
                                                                                       \end{array}
                                                                                     \right).$ This example also shows that (\ref{coerexpaff}) is not sufficient for constant exponential stabilizability of (\ref{SS}).

\item Note that (\ref{bal}) (and so is (\ref{coerexp})) is not necessary for constant uniform exponential stabilization (even if $S_0(t)$ is of isometries and $B$ is positive) as we can see for $A=0$ and $B=\left(
                                                                                                                      \begin{array}{cc}
                                                                                                                        1 & 4 \\
                                                                                                                        0 & 4 \\
                                                                                                                      \end{array}
                                                                                                                    \right)
 $ in $H=R^2.$

\item    If $B\in {\cal L}(H), $ then under the assumptions of Theorem \ref{thm-suffunif}, the system $\dot{y}(t)=Ay(t)+v(t) (\gamma B+\mu B^*)y(t)$ is uniformly exponentially stabilizable by the constant control for all $\gamma,\mu\ge 0$ such that $ (\gamma,\mu)\ne(0,0)$. Moreover,  if $B\ge 0$  and if $\dot{y}=Ay+v(t) (B+B^*)y$ is  exponentially stabilizable with a constant control, then so is both the system  $\dot{y}=Ay+v(t) By$ and $\dot{y}=Ay+v(t) B^*y$. However, the converse is not true as we can see by taking : $A=0$ and $B=\left(
                                                                                                                      \begin{array}{cc}
                                                                                                                        1 & 4 \\
                                                                                                                        0 & 4 \\
                                                                                                                      \end{array}
                                                                                                                    \right)
 $ in $H=R^2.$

 \end{enumerate}
\end{rem}

\subsection{Bilinear systems: Necessary conditions revisited }
In this subsection, we deal with bilinear systems and we will see that, for this class of systems, the observability assumption (\ref{coerexp}) is necessary  for uniform exponential stabilizability  with  constant controls.
 Note that in the case where $B\in {\cal L}(H)$ is self-adjoint and positive, the inequalities (\ref{coer}), (\ref{coerexp}) and (\ref{coerexpaff})
are equivalent to the following one :
\begin{equation}\label{coerexp'}
 \displaystyle\int_0^T \|B^{\frac{1}{2}} S_0(t)y\|^2 dt \ge \delta \|y\|^2,
\;\forall y \in H,
\end{equation}
which means  that the system : $  \dot{\phi}(t)=A\phi(t),$ augmented with the output : $y(t)=B^{\frac{1}{2}}\phi(t),$ is  observable on $[0,T],$ or, equivalently, that the  dual system : $\dot{\xi}(t)=A^*\xi(t)+B^{\frac{1}{2}} u(t)$ is exactly controllable on $[0,T]$  (see \cite{zab95}).

In the next result, we will show that the assumption (\ref{coerexp}) is necessary for uniform exponential stabilization of conservative bilinear systems.

\begin{thm}\label{thm-bil}
Let $A$ generate a linear $C_0-$semigroup  of isometries $S_0(t)$ on $H,$ and let  $ B$ be a linear bounded positive operator.
 Then the assumption (\ref{coerexp}) is necessary  for exponential stabilization of (\ref{SS}) with a constant control.
\end{thm}

\begin{pf}

Applying  Theorem \ref{thm-necaffunif}, we deduce that the estimate (\ref{coerexpaff}) holds. Thus, we may distinguish the following cases :\\
{\it Case 1 : } If $B=B^*\ge 0$, then (\ref{coerexpaff}) is equivalent to (\ref{coerexp}) which gives the claimed result.\\
{\it Case 2 : } Let us return to the case of non self-adjoint control operators. We will apply the first case to the operator $B+B^*.$ Let $v(t)=-\lambda, \; \lambda>0 $ be an exponentially stabilizing control for (\ref{SS}), and let  $M\ge 1,\, \sigma>0 $
be such that the corresponding mild solution  $ S_\lambda(t)z_0$ of (\ref{SS}) satisfies (\ref{exp-est}). \\
Let $T_\lambda(t)$ be the linear semigroup generated by the operator $A-\lambda(B+B^*).$ Since the semigroup $S_\lambda(t)$ is of contractions, and since (\ref{nc0}) holds,  Theorem \ref{thm-suffunif} guarantees the existence of $\tilde{M}, \tilde{\sigma}>0$ such that
\begin{equation}\label{exp-estT}
\|T_\lambda(t)z_0\|\le \tilde{M} e^{-\tilde{\sigma} t} \|z_0\|,\, \forall z_0\in H.
\end{equation}
Since the operator $B+B^*$ is self-adjoint and positive, we deduce from the first case that the operator  $B+B^*$ verifies (\ref{coerexp}) for some $\tilde{T},\tilde{\delta}>0,$ and hence $B$ verifies (\ref{coerexp})
for $T=\tilde{T}$ and $\delta=\displaystyle\frac{\tilde{\delta}}{2}.$
\end{pf}

\begin{rem}
\begin{enumerate}
\item If  $-B$ is positive, then the necessity of (\ref{coerexp}) should be imposed on $-B$.

\item As a consequence of the above result, we retrieve that a linear system with  isometric C$_0-$semigroups  can not be uniformly exponentially stabilized under   compact perturbation of the generator \cite{GUO}.

    \end{enumerate}

\end{rem}

\section{Finite-dimensional systems}

In this subsection, the system (\ref{SS}) is considered in the  Euclidean space $H=R^n $ (with the conventional inner product) and $ A $ is a matrix satisfying the following LMI (linear matrix inequality):
 \begin{equation}\label{Eq-L}
A^TP + PA \le 0,
\end{equation}
for some matrix $P=P^T>0 : $

\begin{thm}\label{thm-Rn1}
1) Let $B$ be Lipschitz,   let $A$ satisfy  (\ref{Eq-L}), and assume that

\begin{equation}\label{coer-n}
 \displaystyle\int_0^T \langle PB e^{tA}y,e^{tA}y\rangle dt \ge \delta \|y\|^2, \;\forall y \in H,
\end{equation}
for some  $ T, \;\delta> 0.$ Then the system (\ref{SS}) is uniformly exponentially stabilizable with a constant control.\\
2) Suppose that there exists $P=P^T > 0$ solution of the following Lyapunov equation:
\begin{equation}\label{lyap-Eq}
 A^TP + PA = 0\cdot
 \end{equation}
  If (\ref{SS}) is constant exponentially stabilizable, then
 $$
  \displaystyle\int_0^T  \|B e^{tA}y\| dt \ge \delta \|y\|,
\;\forall y \in H,
$$
for some $ T, \;\delta> 0.$

\end{thm}

\begin{pf}

1) Since $H$ is of finite dimension, we can replace the inner product $\langle\cdot,\cdot\rangle$ by the one defined by $\langle y,z \rangle_P: = \langle Py,z \rangle$ with corresponding norm $\|\cdot\|_P.$ We have by  (\ref{Eq-L}) that : $\|e^{tA}y\|_P\le \|y\|_P,\, \forall y\in H.$ In other words, $S_0(t)=e^{tA}$ is a contraction semigroup with respect to the new inner product $\langle\cdot,\cdot\rangle_P.$ Then, to apply Theorem \ref{thm-suffunif}, it suffices to observe  that all norm in a finite-dimensional space are equivalent and that  assumption (\ref{coer-n}) means that   (\ref{coerexp}) holds for the  inner product $\langle\cdot,\cdot\rangle_P.$ \\
2) It follows from (\ref{lyap-Eq}) that the semigroup $S(t)=e^{tA}$ is of isometries with respect to the inner product
$\langle\cdot,\cdot\rangle_P. $
 Hence, according to  Theorem \ref{thm-necaffunif}, the estimate (\ref{coerexpaff}) holds in $(H,\langle \cdot,\cdot\rangle_P).$ In the other words, we have:
 $$\displaystyle\int_0^T \langle PB e^{tA}y,Be^{tA}y\rangle dt \ge \delta \langle Py,y\rangle,
\;\forall y \in H,
$$
for some $ T, \;\delta> 0,$ which is equivalent to
 $$
  \displaystyle\int_0^T  \|B e^{tA}y\| dt \ge \delta \|y\|,
\;\forall y \in H,
$$
for some $ T, \;\delta> 0.$
\end{pf}

In the next result, we  study the constant stabilization of a finite-dimensional bilinear system under an algebraic assumption. Let  matrices  $ A $ and $B$ be such that : there exists a non-empty set $\Omega\subset R^n-{\{0\}}$, which complement $\Omega^c$ with the following algebraic property : for each $y\in \Omega, $ there exists $ k\in N$ such that
\begin{equation}\label{lie}
span\{Ay,ad^0_A(B)y,ad^1_A(B)y,..,ad_A^k(B)y\}=R^n,
\end{equation}
where $ ad^0_A(B)=B,ad^1_A(B)= [A,B] = AB-BA, $ and for $k\in I\!\!N, \; ad_A^{k+1}(B) = ad_A(ad_A^k(B)).$\\
The algebraic assumption (\ref{lie}) is equivalent to the following temporal version (see \cite{jurj,sle78}) :

\begin{equation}\label{ad1}
 \langle Be^{tA}y,e^{tA}y\rangle = 0, \; \forall t\ge0
\Longrightarrow y = 0,
\end{equation}
where $\langle\cdot,\cdot \rangle$ is an  inner product in $H=R^n.$

We also have the following lemma that gives link between  (\ref{bal}) and (\ref{coer})  in the context of finite dimensional state spaces.

\begin{lem} \label{lem-*}
Let  $A,  B\in {\cal L}(H)$. If $\dim H<\infty,$  then the
assumptions (\ref{bal}) and (\ref{coer}) are equivalent.
\end{lem}

\begin{pf}  It is clear that (\ref{coer}) $\Rightarrow $ (\ref{bal}).\\
Suppose that (\ref{bal}) holds, and assume by
contradiction that (\ref{coer}) does not hold. Then, for all  $t>0$ and for all integer
 $k\ge 1, $ there exists $z_k\in H$ such
that $\|z_k\|=1$ and
$$\displaystyle\int_0^t
|\left<Be^{sA}z_k,e^{sA}z_k\right>| ds \to 0, \; \mbox{as}\; k\to +\infty. $$

Since $\dim(H)< \infty, $  there exists a convergent subsequence of
$(z_k)$, still denoted by $(z_k)$, and let
$z=\displaystyle\lim_{k\to +\infty} z_k.$\\
For all $t\ge 0,$ we have:\\

$
  | \langle Be^{tA}z_k,e^{tA}z_k\rangle-\langle Be^{tA}z,e^{tA}z\rangle | \le |\langle Be^{tA}(z_k-z),e^{tA}z_k\rangle|+ |\langle Be^{tA}z,e^{tA}(z_k-z)\rangle| \le
  C\|z_k-z\|,\, C>0,
 $\\

from which we can deduce that for all $t\ge 0,$ $$\displaystyle\int_0^t |\left<B
e^{sA}z_k,e^{sA}z_k\right>|dt \to \displaystyle\int_0^t |\left<B
e^{sA}z,e^{sA}z\right>| ds, \,\, \mbox{as}\,\, k\to \infty\cdot$$
We conclude that
$$\displaystyle\int_0^t| \left<Be^{sA}z,e^{sA}z\right>| ds=0,\, \forall t\ge 0. $$
Thus, since the map $t \mapsto |\left<Be^{tA}z,e^{tA}z\right>|$ is continuous, we deduce that
$$ \left<Be^{tA}z,e^{tA}z\right>=0,\; \forall t\ge 0\cdot$$
Then (\ref{bal}) yields $z=0, $ which is a
contradiction, and this achieves the proof.
\end{pf}

\begin{rem}
If $B$ is linear and such that $PB\ge 0$, then (\ref{coer-n})
is equivalent to : $
 \langle PB e^{tA}y,e^{tA}y\rangle=0,\; \forall t\ge 0 \Longrightarrow y=0.
$\\
Recall  that if $B$ is linear,  positive and commutes with  $P$, then $PB$ is also positive.\\
\end{rem}

From the above discussion, we can formulate the following result  which is a consequence of  Theorem \ref{thm-Rn1} and Lemma \ref{lem-*}:

\begin{cor}\label{thm-Rn}
Let $B\in {\cal L}(H).$\\
1) Let   (\ref{Eq-L}) and (\ref{lie}) hold and assume that $PB$  has constant sign. Then the system (\ref{SS}) is constant exponentially stabilizable.\\
2) Suppose that (\ref{lyap-Eq}) holds. Then the condition (\ref{lie}) is necessary  for exponential stabilization of (\ref{SS}) with a constant control.
\end{cor}

\section{Robustness}
In this subsection, we discuss the robustness of
the controller (\ref{cte})  with respect to perturbations of the
parameters $A$ and $B$ of (\ref{SS}). More precisely, we will  exhibit a class of admissible
perturbations of the system (\ref{SS}) that leave its uniform exponential
stability by the control  (\ref{cte}) unaffected.

\subsection{The case of nonlinear control operator}
Let us reconsider the system (\ref{SS}) with a nonlinear control operator $B, $ and let us consider the perturbed system :
\begin{equation}\label{SSA1}
\displaystyle\frac{d z(t)}{d t}= (A+a)z(t) + v(t)Bz(t),\;\;z(0) =
z_0,
\end{equation}
where $ a : H \rightarrow H$ is a linear bounded operator, which represents a perturbation of  (\ref{SS}) on its dynamic $A$. We consider the following set of linear perturbations : ${\cal P}_L=\{a\in{\cal L}(H);\; \langle az,z\rangle\le 0,\;\forall z\in H\}.$ \\
In the sequel, for any operator, $N : H \rightarrow H$, which is Lipschitz and vanishes at $0$, we set  $L_N:=\sup_{y\ne 0}{\frac{\|N (y)\|}{\|y\|}}.$\\

\begin{thm}\label{thmRa1}
Let  assumptions of Theorem \ref{thm-suffunif} hold.
Then  the control (\ref{cte}) uniformly  exponentially stabilizes (\ref{SSA1}) for any perturbation $a\in
{\cal P}_L$ such that  $\|a\| <r:=\displaystyle\frac{-1+\sqrt{1+\frac{\delta}{TL_B}-\frac{\lambda \tilde{K}}{L_B}}}{T},$ where $\tilde{K}$ is the parameter defined in the proof of Theorem \ref{thm-suffunif}.
\end{thm}

\begin{pf} Let $a\in {\cal P}_L.$ Since $a$ is a linear bounded and dissipative
operator, the operator $A+a$ is the infinitesimal generator of a
semigroup of contractions $S_a(t)$ given,  for all  $ y\in
H$, by
$$
S_a(t)y=S_0(t)y+\int_0^tS_0(t-s)aS_a(s)y ds, \, t\ge 0.
$$
Then we have
$$
\langle BS_a(t)y,S_a(t)y\rangle= \langle
BS_0(t)y,S_0(t)y\rangle+\alpha(t)
$$
 where $\alpha(t)$ is a scalar valued function, which is such that
$$
|\alpha(t)|\le TL_B \|a\|\left (2+T \|a\|\right )
\|y\|^2,\; \forall t\in [0,T].
$$
Then, it follows from (\ref{coerexp}) that
\begin{equation}\label{***}
\int_0^T\langle BS_a(t)y,S_a(t)y\rangle dt\ge \delta_a \|y\|^2
\end{equation}
where $\delta_a = \delta-T^2L_B \|a\|(2+T \|a\|).$\\
Applying Theorem \ref{thm-suffunif}, we have exponential stabilization with any constant control (\ref{cte}) such that $0<\lambda < \inf{(\eta_1,\eta_2,\frac{\delta_a}{T\tilde{K}})}$. Let us show that (\ref{cte})  is a common stabilizing control for all perturbed systems (\ref{SSA1}). Remarking that in the proof of Theorem \ref{thm-suffunif} the constants $\eta_i,\, i=1, 2$ do not depend on $\delta$, we can see that it suffices to look for a  $r>0$ for which:
$$
0<\lambda <\frac{\delta_a}{T\tilde{K}},\; \forall a\in {\cal P}_L; \; \|a\|<r
$$
i.e., $
X :=\|a\|<r \; \Rightarrow P(X):=(T^3L_B) X^2+2(T^2L_B) X +\lambda T\tilde{K} - \delta<0.$

By a simple computation, we can check that this implication is satisfied for $r=\frac{-1+\sqrt{1+\frac{\delta}{TL_B}-\frac{\lambda \tilde{K}}{L_B}}}{T},$ and this completes the proof.
\end{pf}
Let us now consider the problem of robustness associated to
perturbations acting, jointly, on the dynamic and the operator of
control.  Consider the perturbed system :

\begin{equation}\label{SSAB1} \displaystyle\frac{d z(t)}{d t}= (A+a)z(t)
+ v(t)(B+b)z(t),\;\;z(0) = z_0,
\end{equation}
 where $a\in{\cal L}(H)$ and $ b : H \rightarrow H$ is a Lipschitz operator such that $b(0)=0,$ so that $0$ remains an equilibrium of (\ref{SSAB1}).\\
Let us introduce the following set of  nonlinear perturbations : ${\cal P}_N =\{n : H\to H; \, n $ is Lipschitz and $ n(0)=0\}.$

\begin{thm}\label{thmRab1}
Let  assumptions of Theorem \ref{thm-suffunif} hold.
Then  the control (\ref{cte}) is robust under any perturbation  $(a, b)\in {\cal P}_L\times {\cal P}_N$ such
  that $\|a\| \le \tilde{r}<r:=\displaystyle\frac{-1+\sqrt{1+\frac{\delta}{TL_B}-\frac{\lambda \tilde{K}}{L_B}}}{T}$ and  $
L_b<-\frac{P(\tilde{r})}{T}, $ where $P(X):=(T^3L_B) X^2+2(T^2L_B) X +\lambda T\tilde{K} - \delta.$
\end{thm}

\begin{pf}
Let $(a, b)\in {\cal P}_L\times {\cal P}_N$, and let $S_a(t)$ denote the semigroup generated by the
operator $A+a$. Since $S_a(t)$ is of contractions, we have from (\ref{***})
$$ \displaystyle\int_0^T\left<(B+b)S_a(t)y,S_a(t)y\right> dt \ge
(\delta_a -TL_b)\|y\|^2. $$
Thus, as in the proof of Theorem \ref{thmRa1}, we should consider the following inequality $P(X)+TL_b<0, \, \forall X\in (0,\tilde{r}),  $ where $P(X):=(T^3L_B) X^2+2(T^2L_B) X +\lambda T\tilde{K} - \delta$. \\
Since $P(\tilde{r})=\min_{X\in (0,\tilde{r})} P(X) <0, $ it suffices to have $P(\tilde{r})++TL_b<0, $ i.e. $L_b<-\frac{P(\tilde{r})}{T}.$
\end{pf}

\subsection{The case of linear control operator}

Here, we reconsider the bilinear case. Let us consider the following perturbed system
\begin{equation}\label{SSA}
\displaystyle\frac{d z(t)}{d t}= Az(t) + v(t)Bz(t)+nz(t),\;\;z(0) =
z_0,
\end{equation}
where $ n : H \rightarrow H$ is a (possibly) nonlinear operator.\\
Let us also consider the perturbed system on $A$ and $B$ :

\begin{equation}\label{SSAB} \displaystyle\frac{d z(t)}{d t}= (A+a)z(t)
+ v(t)(B+b)z(t),\;\;z(0) = z_0,
\end{equation}
 where $a, b : H \rightarrow H$ are (possibly) nonlinear  perturbation operators.

\begin{thm}\label{thmRa}
Let  $A$ generate a semigroup  $S_0(t)$ of contractions on $H,\; B\in
{\cal L}(H)$ and let (\ref{coerexp}) hold. Then,\\
1) the control (\ref{cte}) uniformly exponentially stabilizes the system (\ref{SSA}) for any perturbation $n\in {\cal P}_N$ such that $L_n <\frac{\sigma}{M}, $ where $M, \sigma$ are given by (\ref{exp-est}),\\
2)  the control (\ref{cte}) uniformly exponentially stabilizes (\ref{SSAB}) under any perturbation  $(a, b) \in {\cal P}_N^2$ such
  that :\\ $ L_a+\lambda L_b< \frac{\sigma}{M}.$
 \end{thm}

\begin{pf} 1) Applying Theorem \ref{thm-suffunif}, we deduce that the estimate (\ref{exp-est}) holds for any solution of the nominal system (\ref{SS}), controlled with (\ref{cte}).\\
The mild solution of the perturbed system (\ref{SSA}) is given by
$$
z_n(t)=S_\lambda (t)z_0+\int_0^tS_\lambda(t-s)n z_n(s) ds,
$$
where $S_\lambda (t)$ is the linear semigroup generated by $A-\lambda B. $\\
Since $a$ is Lipschitz and $a(0)=0,$ the Gronwall inequality yields
$$
\|z_n(t)\|\le M e^{-(\sigma -ML_n)t} \|z_0\|,\; \forall t\ge 0.
$$

This gives the claimed result.\\
2) It follows from the first point, by taking $n=a-\lambda b$ and remarking that $L_{a-\lambda b}\le L_a+\lambda L_b$.
\end{pf}

\begin{rem}\label{rem-robus}
\begin{enumerate}
\item The finite dimensional version of Theorem \ref{thmRa}  has been considered in \cite{zab95}, where the nominal system is perturbed with a $C^1-$function $n$ from $H:=R^n $ to $H$.

\item For $0<\lambda<\frac{\delta}{T\tilde{K}},$ one can get a non coupled condition  (separate bounds) on $L_a$ and $L_b$ in the second point of the above theorem, i.e.  one can take $L_a<\frac{\sigma}{2M}$ and $L_b<\frac{\sigma T\tilde{K}}{2M\delta}.$

\end{enumerate}
\end{rem}

\section{Applications}

\begin{ex}
Let us consider the following undamped wave equation
$$
y_{tt}(t)=\Delta y(t)+u(t) y(t)
$$
defined on $\Omega\subset R^n, n\ge 1,$ with Dirichlet boundary conditions.\\
Here $B=\left(
          \begin{array}{cc}
            0 & 0 \\
            I & 0 \\
          \end{array}
        \right)
$ is compact in $X:=H_0^1(\Omega)\times L^2(\Omega)$, so (\ref{coerexpaff})  is impossible (see e.g. \cite{bal1,ouz10}). Then, since $A=\left(
          \begin{array}{cc}
            0 & I \\
            \Delta & 0 \\
          \end{array}
        \right)$ generates a semigroup  of isometries on $H$, the undamped wave equation can not be exponentially stabilizable with none of the controls (\ref{e(t)}), (\ref{sw(t)}) or (\ref{cte}).

\end{ex}

\begin{ex}: Transport equation
\hspace{1cm}\\
Consider the system defined on $H = L^2(0,\infty) $  by the following equation
\begin{equation}
\left \{ \begin{array}{ll}\label{trans}
z_t(\cdot,t)= -z_x(\cdot,t) +v(t)Bz(\cdot,t), & \mbox{in}\;Q\\
&\\
z(\cdot,0) =  z_0,  & \mbox{in}\; \Omega
\end{array}
\right.
\end{equation}
where  $\Omega==(0,\infty),\; Q=(0,\infty)^2$ and $B\in{\cal L}(H).$ Here, we take  $$ Az =
-z_x, \;\forall z\in {\cal
D}(A) = \{ z \in H^1(\Omega ); \;  z(0)= 0 \}.
$$
The operator $ A  $  generates the semigroup of contractions $ S_0(t), t\geq 0 $ defined, for all $z_0\in H,$  by (see e.g. \cite{trig92})
$$ (S_0(t)z_0)(x)=
\left\{\begin{array}{ll}
z_0(x-t)&, \mbox{if} \;\; x > t \\
0 &, \mbox{if } \;\; x \le t
\end{array}
\right.
$$

{\bf Case 1.} $B=g(\cdot) I,$ with $g\in L^\infty(0,\infty).$\\
In the sequel, we take  $g\in L^\infty(0,+\infty)$ such that $ g \ge c>0 $ in $(3,+\infty) $ and  $ g(x)= \left\{
         \begin{array}{ll}
           -1, & 0<x<1 \\
\\
           1, & 1\le x\le 3
         \end{array}
       \right.$\\
We will establish (\ref{coerexp}) for $T=3$. We have

$
  \int_0^3\langle BS_0(t)z_0,S_0(t)z_0 \rangle dt  =  \int_0^3 \int_t^\infty g(x)z_0(x-t)^2 dx dt=
   \int_0^\infty  z_0(x)^2 \int_0^3 g(x+t) dt dx=    \int_0^1z_0(x)^2(\int_0^{1-x} (-1) dt + \int_{1-x}^3 dt) dx +
    \int_1^\infty z_0(x)^2 (\int_0^3 g(x+t) dt) dx \ge     \int_0^1z_0(x)^2[(-1+x)+ (2+x)] dx + 3c \int_1^\infty z_0(x)^2) dx.
$

Thus
$$\int_0^3\langle BS_0(t)z_0,S_0(t)z_0\rangle dt\ge M\|z_0\|^2,\;
M=3\min(1,c),
$$
which implies (\ref{coerexp}).\\
This example shows that (\ref{coerexp}) may hold without the positivity of $B$.\\
Note that the uncontrolled system (\ref{trans}) $ (i.e. \;v(t) =
0) $ is not strongly stable, but it is weakly stable since we have \cite{trig92}
$$
\begin{array}{ccc}
  \langle S_0(t)z_0,y \rangle  &=& \displaystyle
\int_t^{+\infty}z_0(x-t)\overline{y}(x) dx \\
\\
  &\le & \|z_0\|(\displaystyle
\int_t^{+\infty}|y(x)|^2dx)^{\frac{1}{2}}\rightarrow 0, \;\mbox{as}\; t\rightarrow +\infty,
\end{array}
$$
so the constant control (\ref{cte}) enables us to improve the degree of stability
of the system (\ref{trans}).\\

{\bf Case 2.} Let us now give a negative stabilization result regarding the system (\ref{trans}) : consider the linear control operator $$By=\sum_{j\ge1} \frac{1}{j} \langle y, \phi_j\rangle \phi_j,$$
with $\phi_j=\sqrt2\sin(\pi x){\bf 1}_{(0,1)}$, where ${\bf 1}_{(0,1)}$ indicates the characteristic function of $(0,1)$. The operator $B$  is compact, and hence (\ref{coerexp}) is not verified, as we can see by taking $y=\phi_j,\; j\ge 1$ in (\ref{coerexp}). According to Theorem \ref{thm-necaffunif} we deduce that the system (\ref{SS}) can not be  exponentially stabilized by any constant control.
\end{ex}

\begin{ex} Wave equation

Let $\Omega $ denote a bounded open subset
of $ I\!\!R^{n},\, n\ge 1$ with  $C^\infty$ boundary $\partial\Omega$, let
$Q=\Omega\times (0,+\infty) $ and $ \Sigma=\partial \Omega\times  (0,+\infty)$, and consider the following boundary value
problem for the wave equation :\\
\begin{equation}\label{wave}
    \left\{%
\begin{array}{ll}
z_{tt}(x,t)=\Delta z(x,t)+v(t) G(z_t)(x,t), & in\ Q \\
    \\
    z(\xi,t)=0, & on\ \Sigma\\
\end{array}%
\right.
\end{equation}
where $G : L^2(\Omega) \rightarrow L^2(\Omega)$ is a Lipschitz operator such that $G(0)=0,$ and $v(t)$ is the multiplicative control $v(t)$.
 When $v(t) = 1,$ this problem has been investigated extensively under
the assumption that the norm of the indefinite damping operator is small enough
 (see \cite{ben,cheG,dao,frei1,frei2,las2} and the references therein).\\
   The system (\ref{wave}) has the form of equation (\ref{SS}) if we set
 $$ A=\left(%
\begin{array}{cc}
  0 & I \\
  \Delta & 0 \\
\end{array}%
\right)\ and \ B\left(%
\begin{array}{c}
   y\\
   z \\
\end{array}
\right)=\left(%
\begin{array}{c}
   0\\
   G(z) \\
\end{array}
\right).$$
Here, we consider the Hilbert state space $H=H^{1}_{0}(\Omega)\times L^{2}(\Omega)$ with the inner product
$\left<(y_1,z_1),(y_2,z_2)\right>=\left <y_1,y_2\right>_{H_0^1(\Omega)}+\left <z_1,z_2\right>_{L^2(\Omega )}$.
Thus, the operator $A$ with domain $D(A)=(H^{2}(\Omega)\cap H^{1}_{0}(\Omega))\times H^{1}_{0}(\Omega)$ is skew-adjoint,
so the corresponding semigroup  $S_0(t)$ is of isometries. Since $G$ is supposed Lipschitz,  the operator $B$ is Lipschitz from $H$ to $H$.
Then, according to Theorem \ref{thm-suffunif}, the system (\ref{wave}) is uniformly constant exponentially stabilizable if the assumption (\ref{coerexp}) holds,
which is equivalent to
\begin{equation}\label{coerG-wave}
\begin{array}{cc}
  \int_\Omega \big (|\varphi_t(x,0)|^2+ |\nabla \varphi(x,0)|^2 \big) dx \le &  \\
    C \int_0^T \langle G(\varphi_t)(x,t),\varphi_t(x,t)\rangle dx dt,&
\end{array}
\end{equation}
for some $ C>0, $ where $\varphi$ is the solution of the uncontrolled wave equation
$$    \left\{%
\begin{array}{ll}
    \varphi_{tt} (x,t)=\Delta \varphi(x,t), & on\ Q \\
    \\
  \varphi(0)=\varphi_1,\; \varphi_t(0)=\varphi_2\\
  \\
    \varphi=0, & on\ \partial \Omega\times  (0,+\infty)\\
\end{array}%
\right.
$$
Note that $$
S_0(t)z_0=(\varphi,\varphi_t),\; z_0=(\varphi_0,\varphi_1)\in H.
$$
Moreover, a necessary condition for constant uniform exponential stabilization of system (\ref{wave}) is given by :
$$
\big ( \int_\Omega |\varphi_t(x,0)|^2+ |\nabla \varphi(x,0)|^2 dx \big)^{\frac{1}{2}}  \le
$$
$$
C \int_0^T \|G(\varphi_t)(x,t)\|  dt, \; (C>0).$$
Let us now give some examples for the damping $G.$\\
{Case 1 : } Let $\omega$ be a non-empty open subset of $\Omega$,  let $\phi$ be an unitary vector of $L^2(\Omega), $ and let us consider the nonlinear operator:
$$
G(y)=\big (|\int_\omega y(x)\phi(x) dx|\phi+2y\big ) {\bf 1}_\omega.
$$
We have \\

$\langle BS_0(t)y_0,S_0(t)y_0\rangle=|\langle \varphi_t,\phi|_\omega\rangle| \langle \varphi_t,\phi|_\omega\rangle+2\|\varphi_t|_\omega\|^2.$

Then it follows from Cauchy-Schwartz inequality that $$
\langle BS_0(t)y_0,S_0(t)y_0\rangle\ge \|\varphi_t|_\omega\|^2
$$
Thus,  the estimate (\ref{coerG-wave}) holds if the following one is satisfied :
\begin{equation}\label{coer-wave}
\int_\Omega (|\varphi_t(x,0)|^2+ |\nabla \varphi(x,0)|^2) dx \le C \int_0^T \int_\omega |\varphi_t(x,t)|^2 dx dt.
\end{equation}
In this case, we can conclude the exponential stability by using Theorem \ref{thm-suffunif}.\\
It is well known that (\ref{coer-wave}) holds under a geometrical assumption on $\omega$.
  As an interpretation of the inequality (\ref{coer-wave}), we have the following "geometric control property" :  every ray of geometric
optics propagating in $\Omega$  and being reflected on its boundary, enters the control region $\omega$ in a time less
than $T$ (see \cite{bar,lag}).\\
{\bf Case 2 : } Here, we will consider a linear  damping operator. Let $(\phi_j)_{j\ge 1}$  be an orthonormal basis of $L^2(\Omega)$ formed with eigenfunctions  of the Laplacian operator associated to Dirichlet boundary conditions, and let $\lambda_j,\, j\ge 1$ denote the
corresponding eigenvalues. Let us consider the linear operator :
$$
G(y)=\big (\int_\Omega y(x)\phi_2(x) dx \big )\phi_1+2y,
$$
To establish (\ref{coerexp}), let $y=(y_1,y_2)\in H$ with $y_1=\displaystyle\sum_{j=1}^{\infty}
\alpha_j\phi_j$ and $y_2=\displaystyle\sum_{j=1}^{\infty}
\lambda_j^{\frac{1}{2}} \beta_j\phi_j, \;$ where $ \; (\alpha_j,\beta_j) \in I\!\!R^2, \;  j\ge 1\cdot$\\
We have
$\|y\|^2=\displaystyle\sum_{j=1}^{\infty}\lambda_j(\alpha^2_j +
\beta^2_j).$ Separation of variables yields
$$
S_0(s)y=\displaystyle\sum_{j=1}^{\infty}\left(%
\begin{array}{c}
  \alpha_j\cos(\lambda_j^{\frac{1}{2}} s) + \beta_j\sin(\lambda_j^{\frac{1}{2}} s) \\
  \\
 -\alpha_j\lambda_j^{\frac{1}{2}}\sin(\lambda_j^{\frac{1}{2}} s) +
 \beta_j\lambda_j^{\frac{1}{2}}\cos(\lambda_j^{\frac{1}{2}} s)  \\
\end{array}%
\right)\phi_j\cdot
$$
Then we have
$$
  \langle BS_0(s)y,S_0(s)y\rangle \ge \|\varphi_t\|^2=
  $$
  $$\displaystyle\sum_{j=1}^{\infty}\lambda_j\left\{\alpha_j^2\sin^2(j\pi
s) + \beta_j^2\cos^2(j\pi s) - \sin(2j\pi s)
\alpha_j\beta_j\right\}\cdot
$$
It follows that
$$
\displaystyle\int_0^2\langle BS_0(s)y,S_0(s)y\rangle ds
\ge\displaystyle\sum_{j=1}^{\infty}\lambda_j(\alpha^2_j +
\beta^2_j),
$$
so the assumption (\ref{coerexp}) holds for $T=2 $ and  $ \delta= 1$. It follows  that the system (\ref{wave}) is uniformly  exponentially stabilizable by means of a constant control.\\
Note that, $B$ is not self-adjoint, so one can not apply the linear version given in \cite{cur06,har90,las03}.\\
Let us now examine the robustness of the constant  control.  Consider the following perturbed system
\begin{equation}\label{waverobut}
    \left\{%
\begin{array}{ll}
    z_{tt}(x,t)=\Delta z(x,t)+p(x)|z(x,t)|+&\\
        v(t)q(x)z_t(x,t)+ G(z_t)(x,t), & on\ Q \\
    \\
    z=0, & on\ \Sigma\\
\end{array}%
\right.
\end{equation}
where $ p, q\in L^\infty(\Omega)$.\\
From Theorem \ref{thmRa} and Remark \ref{rem-robus}, we have exponential stability of (\ref{waverobut}) with the constant control (\ref{cte})  provided that $\|p\|_{L^\infty(\Omega)}<\frac{\sigma}{2M}=-\frac{\sqrt \gamma \ln(\gamma)}{4}$ and $\|q\|_{L^\infty(\Omega)}<\frac{\sigma T\tilde{K}}{2M\delta}=-\sqrt \gamma \ln(\gamma),$ where $\gamma$ is given, according to in the proof of Theorem \ref{thm-suffunif}, by : $\gamma=8\lambda^2-2\lambda+1\in (0,1), \; $ for all $ 0<\lambda<\frac{1}{4}.$\\

\end{ex}

\section{Conclusion}
This paper studied  exponential stabilization of control-affine systems using constant controls.
The main assumptions of necessity and sufficiency are formulated in terms of observability like estimates. As an important subclass of control-affine systems,  the class of bilinear systems with  bounded control operators. However, the modelization can give rise to the unboundedness aspect of the operator of control. This is the case of  systems  with boundary or pointwise controls. The stabilization problem of unbounded bilinear systems  has been considered in \cite{ber10,elay12} using  nonlinear feedback controls. Now, a  natural question is to ask whether an unbounded bilinear  system can be stabilized by means of a constant control.

\end{document}